\documentclass[12pt, a4paper]{article}
\usepackage[english]{babel}
\usepackage{a4wide}
\usepackage{graphicx}
\usepackage{amssymb,amsmath,amsthm}
\usepackage{amsfonts}
\usepackage{array}
\newtheorem{theorem}{Theorem}
\newtheorem{lemma}{Lemma}

\title{On rate of convergence estimates for~homogeneous discrete-time nonlinear Markov chains}
\author{Aleksandr Shchegolev\footnote{National Research University Higher School of Economics,
Moscow, Russia,\linebreak email: ashchegolev@hse.ru. Supported by grant Russian Foundation for Basic Research 20-01-00575.}}
\date{March, 2021}

\begin{document}
\maketitle

\abstract{The paper studies an improved estimate for the rate of convergence for nonlinear homogeneous discrete-time Markov chains. These processes are nonlinear in terms of the distribution law. Hence, the transition kernels are dependent on the current probability distributions of the process apart from being dependent on the current state. Such processes often act as limits for large-scale systems of dependent Markov chains with interaction. The paper generalizes the convergence results by taking the estimate over two steps. Such an approach keeps the existence and uniqueness results under assumptions that are analogical to the one-step result. It is shown that such an approach may lead to a better rate of convergence. Several examples provided illustrating the fact that the suggested estimate may have a better rate of convergence than the original one. Also, it is shown that the new estimate may even be applicable in some cases when the conditions of the result on one step cannot guarantee any convergence. Finally, these examples depict that the original conditions may not be an obstacle for the convergence of nonlinear Markov chains.}


\section{Introduction}
\label{S:1}

In earlier papers ergodic properties of ordinary Markov chains were studied by lots of authors; we only mention a few, A.A.~Markov, A.N.~Kolmogorov, W.~Doeblin, J.L.~Doob, R.L.~Dobrushin. There exist several extensions to the theory of Markov processes connected with the dependence on the distribution law of the process itself. O.~Onicescu and G.~Mihoc \cite{Onicescu:1935} introduced ``Chains with complete connections'' -- the processes, which depend on the conditional distribution of the previous step. Ergodic properties of such processes were studied by the school of Romanian mathematicians, the main results in that field were summarized in the monograph by M.~Iosifescu and S.~Grigorescu \cite{Iosifescu:1990}. Such generalization is also known and studied in the field of symbolic dynamics under the name of $g$-measures \cite{Keane:1972}. Another extension is connected with the class of processes, introduced by H.P.~McKean \cite{McKean:1966} in 1966. This kind of processes is also studied by lots of researchers including A.-S.~Sznitman \cite{Sznitman:1991} and V.N.~Kolokoltsov \cite{Kolokoltsov:2010}. The case we are considering is included in this extension and is called nonlinear Markov chains. Here, nonlinearity means the dependence of the transition functions of the process on its state and the distribution law at the current moment. Unlike the ``chains with complete connections'' \cite{Onicescu:1935}, these processes assume that the process depends on the unconditional distribution of the process. The results of studying the ergodic properties of homogeneous nonlinear Markov chains in discrete time were obtained by O.A.~Butkovsky \cite{Butkovsky:2014}. In particular, he showed that in the nonlinear case the theory of ordinary Markov chains is insufficient, and an additional condition for the convergence was established. Such processes are interesting as they often occur as limits for large systems of dependent Markov chains with interaction.

In this paper, some generalization of the existing results on convergence for homogeneous nonlinear Markov chains with discrete-time \cite{Butkovsky:2014} is made, using the estimate over several steps. In the first part, a generalization of the estimate for ergodic convergence is obtained, which depends on the transition probabilities over two steps. The second part presents examples of homogeneous nonlinear Markov chains for which the convergence conditions for the new result are satisfied, while the existing result for one step is either inapplicable or leads to slower convergence.

\section{Problem setting and main result}
\label{S:2}
Convergence and uniform ergodicity conditions for discrete irreducible aperiodic homogeneous Markov chains with a finite state space shown in various classical sources, and similar results also exist for more general Markov chains (for example, \cite{Veretennikov:2017}).

Let $(E, \mathcal{E})$ be a measurable space and $\mathcal{P}(E)$ be a set of probability measures defined on this space. Then the process $\left(X_n^\mu\right)_{n\in\mathbb{Z}_{+}}$ is a nonlinear Markov chain with the space state $(E, \mathcal{E})$, initial distribution $\mu = \text{Law}\left(X_0^\mu\right)$, $\mu \in \mathcal{P}(E)$ and transition probabilities $P_{\mu_n}(x, B) = \mathbb{P}_{\mu_n}\left(X_{n+1}^\mu \in B | X_n^\mu = x\right)$, where $x \in E$, $B \in \mathcal{E}$, $n \in \mathbb{Z}_{+}$ and $\mu_n := \text{Law}(X_n^\mu)$. Thus, the transition kernel is dependent not only on the state of the process in the moment $n$, but also on the distribution in that moment.

Assume $\mu, \nu \in \mathcal{P}(E)$, then the total variation distance between two probability measures may be defined as follows:
\vspace{-0.5em}
\begin{equation*}
\|\mu - \nu\|_{TV} = 2\sup_{A\in \mathcal{E}}|\mu(A) - \nu(A)| = \int_E |\mu(dx) - \nu(dx)|.
\vspace{-0.25em}\end{equation*}

According to the results \cite{Butkovsky:2014}, nonlinear Markov chain is uniformly ergodic and the existence and uniqueness of the invariant measure $\pi$ is guaranteed if the chain satisfies the following conditions:
\begin{align}
\sup_{\mu,\nu\in\mathcal{P}(E)}\|P_\mu(x,\cdot) - P_\nu(y,\cdot)\|_{TV} \le 2(1-\alpha),
\label{NonlinearDobrushinAnalog}
\end{align}
where $0 < \alpha < 1, \quad x,y\in E$;
\begin{align}
\|P_\mu(x,\cdot) - P_\nu(x,\cdot)\|_{TV} \le \lambda \|\mu - \nu\|_{TV},
\label{ButkovskyLambda}
\end{align}
where $\lambda\in [0,\alpha],\quad x\in E,\quad \mu,\nu \in \mathcal{P}(E)$.

Then the convergence is exponential, if $\lambda < \alpha$ and 
\vspace{-0.5em}
\begin{equation*}
\|\mu_n - \pi\|_{TV} \le 2(1-(\alpha - \lambda))^n,\quad n\in\mathbb{Z}_+,
\vspace{-0.25em}\end{equation*}
and in case $\lambda = \alpha$ we have linear convergence
\vspace{-0.5em}
\begin{equation*}
\|\mu_n - \pi\|_{TV} \le \frac{2}{\lambda n},\quad n\in\mathbb{Z}_+,
\vspace{-0.25em}\end{equation*}
while for the case $\lambda > \alpha$ in the paper \cite{Butkovsky:2014} there were some counterexamples showing that there might be either no invariant measures or existence of infinite number of measures.

The purpose of this paper is to show that these counterexamples belong to a certain limited class of nonlinear chains, and for other nonlinear Markov chains, the condition $\lambda>\alpha$ is not an obstacle to exponential convergence.

Let us generalize the result of \cite{Butkovsky:2014} using transition kernels over two steps.

\begin{theorem}[Existence and uniqueness of the invariant measure]
Let the process $X$ have a transition matrix over two steps $Q_\mu(x, A): = P_\mu(X_2\in A | X_0 = x)$ and satisfies the following conditions:
\begin{align}
\sup_{\mu,\nu\in\mathcal{P}(E)}\|Q_{\mu}(x,\cdot) - Q_{\nu}(y,\cdot)\|_{TV} \le 2(1-\alpha_2),
\label{NonlinearDobrushinAnalog2}
\end{align}
where $0 < \alpha_2 < 1, \quad x,y\in E$,
\begin{align}
\|Q_{\mu}(x,\cdot) - Q_{\nu}(x,\cdot)\|_{TV} \le \lambda_2 \|\mu - \nu\|_{TV},
\label{Lambda2}
\end{align}
where $\lambda_2 \in [0,\alpha_2],\quad x\in E,\quad \mu,\nu \in \mathcal{P}(E)$,
\begin{align}
\|P_{\mu}(x,\cdot) - P_{\nu}(x,\cdot)\|_{TV} \le \lambda_1\|\mu - \nu\|_{TV},\,\,\lambda_1 < \infty.
\label{Lambda1}
\end{align}
Then for the process $X$ there exists a unique invariant measure $\pi$ and for any probability measure $\mu \in \mathcal{P}(E)$ the following convergence is true
{\small
\begin{flalign}\label{pi_lambda<alpha}
\|\mu_n - \pi\|_{TV} \le \|\mu_0 - \pi\|_{TV}(1-\alpha_2  + \lambda_2)^{[n/2]}((1 + \lambda_1) \mathbf{1}(n \, \text{is odd})\vee 1),
\end{flalign}}
and in case $\lambda_2 = \alpha_2$
\begin{align}
\|\mu_n - \pi\|_{TV} \le \frac{\|\mu_0 - \pi\|_{TV}}{1+ \frac{\lambda_2 n}{2}\|\mu_0 - \pi\|_{TV}}((1 + \lambda_1) \mathbf{1}(n \, \text{is odd})\vee 1).
\label{pi_lambda=alpha}
\end{align}
\label{th1}
\end{theorem}
\vspace{-0.5em}
To prove this theorem, we need an auxiliary theorem on the convergence of any two initial probability measures for a given process.

\begin{theorem}
Let the process $X$ have a transition probability matrix over two steps $Q_\mu(x,B)$ and satisfies conditions \eqref{NonlinearDobrushinAnalog2} and \eqref{Lambda2} of the theorem \ref{th1}. Then for any two probability measures $\mu, \nu \in \mathcal{P}(E)$ the following convergence is true:
{\small
\begin{flalign}\label{lambda<alpha}
\|\mu_n - \nu_n\|_{TV} 
 \le \|\mu_0 - \nu_0\|_{TV}(1-\alpha_2  + \lambda_2)^{[n/2]}((1 + \lambda_1) \mathbf{1}(n \, \text{is odd})\vee 1),
\end{flalign}}
and in case $\lambda_2 = \alpha_2$
\begin{align}
\|\mu_n - \nu_n\|_{TV} \le \frac{\|\mu_0 - \nu_0\|_{TV}}{1+ \frac{\lambda_2 n}{2}\|\mu_0 - \nu_0\|_{TV}}((1 + \lambda_1) \mathbf{1}(n \, \text{is odd})\vee 1).
\label{lambda=alpha}
\end{align}
\label{th2}
\end{theorem}

The proofs of these theorems are similar to \cite{Butkovsky:2014} and using a transition kernel over two steps. The complete proof is presented in order to correct some minor computational inaccuracies of the \cite{Butkovsky:2014} proof, which, however, did not affect the final result.

Let us prove the theorem \ref{th2}.
\begin{proof}
Let $P: E \times \mathcal{E} \rightarrow [0,1]$ be a transition kernel, $\varphi: E\rightarrow \mathbb{R}$ be a measurable function and $\mu \in \mathcal{P}(E)$ be a probability measure; denote $\mu P:= \int_EP(x,dt)\mu(dx)$; in case when $P$ is dependent on the measure $\mu$ we have: $\mu_1(\mu):=\mu P_\mu:= \int_EP_\mu(x,dt)\mu(dx)$, then two-step transition kernel
\vspace{-0.5em}
\begin{equation*}
Q_\mu(x,dy) = \int P_\mu(x,dx_1)P_{\mu_1(\mu)}(x_1,dy).
\vspace{-0.25em}\end{equation*}

Consider the distance between measures in the total variation metric after applying the transition kernel over two steps.

For any probability measures $\mu, \nu \in \mathcal{P}(E)$, we denote $$d\eta = ((d\mu / d\nu) \wedge 1)d\nu$$ and apply the triangular inequality, as a result we get
\begin{align*}
\|\mu Q_{\mu}- \nu  Q_{\nu}\|_{TV} =
\|(\eta + (\mu - \eta)) Q_{\mu}- (\eta - (\nu - \eta))  Q_{\nu}\|_{TV} =\\
= \int_E |\eta Q_{\mu}(dx) + (\mu - \eta)Q_{\mu}(dx)- \eta Q_{\nu}(dx) - (\nu - \eta)Q_{\nu}(dx)| \le\\
\le \int_E |\eta Q_{\mu}(dx) - \eta Q_{\nu}(dx)|+ \int_E |(\mu - \eta)Q_{\mu}(dx) - (\nu - \eta)Q_{\nu}(dx)| \le \\
\le \|\eta Q_{\mu} - \eta Q_{\nu}\|_{TV} + \|(\mu - \eta)Q_{\mu} - (\nu - \eta)Q_{\nu}\|_{TV}.
\end{align*}

Consider the first term, applying the Jensen's inequality and \eqref{Lambda2} to it, and also using the following fact: $\eta(E) = 1 - \|\mu-\nu\|_{TV}/2$. We get
{\small
\begin{align*}
\|\eta Q_{\mu} - \eta Q_{\nu}\|_{TV} &= \int_E \left|\int_EQ_{\mu}(x,dy)\eta(dx) - \int_EQ_{\nu}(x,dy)\eta(dx)\right| \le \\
&\le \int_E\int_E|Q_{\mu}(x,dy) - Q_{\nu}(x,dy)|\eta(dx) \le \\ 
&\le \lambda_2 \|\mu-\nu\|_{TV}\left(1 - \frac{1}{2}\|\mu-\nu\|_{TV}\right).
\end{align*}
}
Then the second term
{\small
\begin{align*}
\|(\mu - \eta)Q_{\mu} - (\nu - \eta)Q_{\nu}\|_{TV} = \int_E |(\mu - \eta)Q_{\mu}(dx) - (\nu - \eta)Q_{\nu}(dx)| =\\
=
\int_E \left|\int_EQ_{\mu}(x,dy)(\mu - \eta)(dx) - \int_EQ_{\nu}(x',dy)(\nu - \eta)(dx')\right|.
\end{align*}
}

Recall that $\mu_n = \text{Law}(X_n^\mu)$, $\nu_n = \text{Law}(X_n^\nu)$, we denote $p_0 = \|\mu_0 - \nu_0\|_{TV}/2$, assuming $p_0>0$ (if $p_0=0$, then $p_2=0$, etc.).

Let us estimate the expression $\|\mu_2 - \nu_2\|_{TV}$ from above.
{\small
\begin{align*}
\|\mu_2 - \nu_2\|_{TV} \le \lambda_2 \|\mu_0 -  \nu_0\|_{TV} \left(1- \frac12 \|\mu_0 - \nu_0\|_{TV}\right)
\\
+ p_0 \int \left|\int Q_{\mu}(x, dy)\frac{(\mu_0-\eta_0)(dx)} {p_0}- \int Q_{\nu}(x', dy)\frac{(\nu_0-\eta_0)(dx')} {p_0}\right|
\\
=  2p_0 \lambda_2 (1-p_0)  + p_0\int
\left|\iint (Q_{\mu}(x, dy) - Q_{\nu}(x', dy)) \frac{(\mu_0-\eta_0)(dx)} {p_0}  \frac{(\nu_0-\eta_0)(dx')} {p_0}\right|
\\
\le \,  2p_0 \lambda_2 (1-p_0)  + p_0\iiint \left| Q_{\mu}(x, dy)- Q_{\nu}(x', dy)\right|\frac{(\mu_0-\eta_0)(dx)} {p_0}\frac{(\nu_0-\eta_0)(dx')} {p_0}
\\
\le 2p_0 \lambda_2 (1-p_0)  +
2(1-\alpha_2)p_0 \iint \frac{(\mu_0-\eta_0)(dx)} {p_0} \frac{(\nu_0-\eta_0)(dx')} {p_0}
\\
= 2p_0 \lambda_2 (1-p_0)  +
2p_0(1-\alpha_2) =
2p_0 (\lambda_2 - \lambda_2p_0 + 1 - \alpha_2).
\end{align*}}
If $\lambda_2 < \alpha_2$, we obtain
\vspace{-0.5em}
\begin{equation*}
\|\mu_2 - \nu_2\|_{TV} \le \|\mu_0-\nu_0\|_{TV}(1 - \alpha_2 + \lambda_2),
\vspace{-0.25em}\end{equation*}
while in case $\lambda_2 = \alpha_2$ we get
\vspace{-0.5em}
\begin{equation*}
\|\mu_2 - \nu_2\|_{TV} \le 2p_0(1 - \lambda_2p_0),
\vspace{-0.25em}\end{equation*}
or
\vspace{-0.5em}
\begin{equation*}
p_2 \le p_0(1 - \lambda_2p_0).
\vspace{-0.25em}\end{equation*}

For the case $2n+1$ we may obtain the following result:
{\small
\begin{align*}
\|\mu_{2n+1} - \nu_{2n+1}\|_{TV} \le \lambda_1 \|\mu_{2n} -  \nu_{2n}\|_{TV} \left(1- \frac12 \|\mu_{2n} - \nu_{2n}\|_{TV}\right)
\\
+ p_{2n} \int \left|\int P_{\mu}(x, dy)\frac{(\mu_{2n}-\eta_{2n})(dx)} {p_{2n}}- \int P_{\nu}(x', dy)\frac{(\nu_{2n}-\eta_{2n})(dx')} {p_{2n}}\right|
\\
=  2p_{2n} \lambda_1 (1-p_{2n})  + p_{2n}\int
\left|\iint (P_{\mu}(x, dy) - P_{\nu}(x', dy)) \frac{(\mu_{2n}-\eta_{2n})(dx)} {p_{2n}}  \frac{(\nu_{2n}-\eta_{2n})(dx')} {p_{2n}}\right|
\\
\le \,  2p_{2n} \lambda_1 (1-p_{2n})  + p_{2n}\iiint \left| P_{\mu}(x, dy)- P_{\nu}(x', dy)\right|\frac{(\mu_{2n}-\eta_{2n})(dx)} {p_{2n}}\frac{(\nu_{2n}-\eta_{2n})(dx')} {p_{2n}}
\\
\le 2p_{2n} (\lambda_1 (1-p_{2n})  + 1) = (1 + \lambda_1)\|\mu_{2n} -  \nu_{2n}\|_{TV}.
\end{align*}}
Hence,
\begin{align}
\|\mu_{2n+1} - \nu_{2n+1}\|_{TV} \le (1 + \lambda_1)\|\mu_{2n} -  \nu_{2n}\|_{TV}.
\label{OneStep1lambda}
\end{align}

Thus, iterating the estimate for $\lambda_2 < \alpha_2$, by induction we obtain
\vspace{-0.5em}
\begin{equation*}
\|\mu_n - \nu_n\|_{TV} \le  \|\mu_0 - \nu_0\|_{TV}(1-\alpha_2  + \lambda_2)^{[n/2]}((1 + \lambda_1) \mathbf{1}(n \, \text{is odd})\vee 1).
\vspace{-0.25em}\end{equation*}

For the case $\alpha_2 = \lambda_2$ we apply the following lemma.
\begin{lemma}
\label{lemma}
Let $a_0, a_1, \dots$ be some sequence of positive numbers. Assume that $0 < a_0 \le 1$ and the following estimate is true
\vspace{-0.5em}
\begin{equation*}
a_{n+1} \le a_n(1-\psi(a_n)),\quad n\in\mathbb{Z}_+,
\vspace{-0.25em}\end{equation*}
where $\psi: [0,\infty)\rightarrow [0,1]$ -- continuous non-decreasing function with $\psi(0) = 0$ and $\psi(x) > 0$ as $x>0$. Then
\vspace{-0.5em}
\begin{equation*}
a_n \le g^{-1}(n)
\vspace{-0.25em}\end{equation*}
for all $n\in\mathbb{Z}_+$, where
\vspace{-0.5em}
\begin{equation*}
g(x) = \int_x^{a_0}\frac{dt}{t\psi(t)},\quad 0 < x \le 1.
\vspace{-0.25em}\end{equation*}
\end{lemma}

This lemma in a slightly different version and its proof are given in \cite{Butkovsky:Thesis}. Since we use a slightly modified version with a different upper limit in the integral, the proof is presented for the convenience of the reader, even though it coincides with the original source.

\begin{proof}
Note that the function $g^{-1}$ exists since $g$ is unbounded, non-negative, and strictly decreasing. Then it follows from the non-negativity of $\psi$ that $a_{n+1}\le a_n$ for any $ n\in\mathbb{N}$. Then there is $s\in[a_{n+1}, a_n]$ such that
$$g(a_{n+1})-g(a_{n}) = g'(s)(a_{n+1}-a_n) = -\frac{a_{n+1}- a_n}{s\psi(s)} \ge \frac{a_n\psi(a_n)}{s\psi(s)} \ge 1.$$
Thus, $g(a_n) \ge n$ and $a_n \le g^{-1}(n)$.
\end{proof}

Applying lemma \ref{lemma} for $a_{n+2}$ and $\psi(t) = \lambda_2 t$ we obtain
\vspace{-0.5em}
\begin{equation*}
p_2 \le g^{-1}(n) = \frac{1}{\lambda_2 n + \frac{1}{p_0}} = \frac{p_0}{1+p_0\lambda_2n}.
\vspace{-0.25em}\end{equation*}

Summarizing the result for odd $n$ using \eqref{OneStep1lambda}, we have:
\vspace{-0.5em}
\begin{equation*}
\|\mu_n - \nu_n\|_{TV} \le \frac{\|\mu_0 - \nu_0\|_{TV}}{1+ \frac{\lambda_2 n}{2}\|\mu_0 - \nu_0\|_{TV}}((1 + \lambda_1) \mathbf{1}(n \, \text{is odd})\vee 1).
\vspace{-0.25em}\end{equation*}
\end{proof}
\vspace{-0.5em}

Next we proceed to to the proof of theorem \ref{th1}.
\begin{proof}
Consider a sequence of probability measures $\left(\mu_n\right)_{n\in\mathbb{N}}$. According to the theorem \ref{th2}, by virtue of \eqref{lambda<alpha} and \eqref{lambda=alpha}, for any $m,n \in \mathbb{N}$
\vspace{-0.5em}
\begin{equation*}
\|\mu_n - \mu_{n+m}\|_{TV} \le \frac{\|\mu_0 - \mu_m\|_{TV}}{1+ \frac{\lambda_2 n}{2}\|\mu_0 - \mu_m\|_{TV}}((1 + \lambda_1) \mathbf{1}(n \, \text{is odd})\vee 1).
\vspace{-0.25em}\end{equation*}
Then $\left(\mu_n\right)_{n\in\mathbb{N}}$ is a Cauchy sequence in the complete metric space $(\mathcal{P}(E), \|\cdot\|_{TV})$ and we may find $\pi\in\mathcal{P}(E)$, such that $\lim_{n\rightarrow \infty}\|\mu_n-\pi\|_{TV} = 0$.

Let us show that the limit measure $\pi$ is invariant. To show this we use the triangular inequality and the condition \eqref{OneStep1lambda} for $ n\rightarrow\infty$:
\vspace{-0.5em}
\begin{equation*}
\|\pi P_\pi - \mu_{n+1}\|_{TV} = \|\pi P_\pi - \mu_{n}P_{\mu_{n}}\|_{TV} \le (1+\lambda_1)\|\pi-\mu_n\|_{TV}\rightarrow 0,
\vspace{-0.25em}\end{equation*}
while $\mu_{n+1} \rightarrow \pi$, we obtain
\vspace{-0.5em}
\begin{equation*}
\|\pi P_\pi - \pi\|_{TV} \le \|\pi P_\pi - \mu_{n+1}\|_{TV} + \|\mu_{n+1} - \pi\|_{TV} \rightarrow 0.
\vspace{-0.25em}\end{equation*}
Hence we get $\pi = \pi P_\pi$.

To prove the uniqueness of the invariant measure $\pi$, assume that $\nu\in\mathcal{P}(E)$ is such that $\nu \ne \pi$ and $\nu = \nu P_\nu$, then iteratively applying the results \eqref{lambda<alpha} and \eqref{lambda=alpha} we get a contradiction
\vspace{-0.5em}
\begin{equation*}
\|\nu-\pi\|_{TV} = \|\nu Q_\nu-\pi Q_\pi\|_{TV} < \|\nu-\pi\|_{TV}.
\vspace{-0.25em}\end{equation*}
Thus, the process $(X_n^{\mu})_{n\in\mathbb{Z}}$ has a unique invariant measure $\pi$.
\end{proof}

In section \ref{S:3} we will illustrate that the obtained estimate can be better than the estimate obtained through one-step transition.

\section{Examples}
\label{S:3}
\subsection{Example of a nonlinear Markov chain, $\alpha > 0, \lambda > \alpha$}
Consider a homogeneous nonlinear Markov chain in discrete time $X_n^\mu$ with state space $(E, \mathcal{E}) = \left(\{1, 2, 3, 4\}, 2^{\{1, 2, 3, 4\}}\right)$, the initial distribution $\mu_0$ and the transition probability matrix $P_{\mu_0}(i, j)$, defined as follows:
\vspace{-0.5em}
\begin{equation*}
\mu_0 = \begin{pmatrix}\mu(\{1\}) & \mu(\{2\})& \mu(\{3\})& \mu(\{4\})\end{pmatrix}
\vspace{-0.25em}\end{equation*}
{\small
\vspace{-0.5em}
\begin{equation*}
P_{\mu_0}(i, j) =
\begin{pmatrix}
0.001 + \gamma\mu(\{1\}) & 0.001 + \gamma\mu(\{1\}) & 0.499 - \gamma\mu(\{1\}) & 0.499 - \gamma\mu(\{1\})\\
0.499 & 0.499 & 0.001 & 0.001 \\
0.499 & 0.001 & 0.499 & 0.001 \\
0.001 & 0.499 & 0.001 & 0.499
\end{pmatrix},
\vspace{-0.25em}\end{equation*}
}
where $0 < \gamma < 0.25$.

We can notice that for a given process the conditions \eqref{NonlinearDobrushinAnalog} and \eqref{ButkovskyLambda} guarantee convergence to an invariant measure only in the case of $\gamma \leqslant 0.004$, since $\alpha = 0.004$ and $\lambda = \gamma$.

Let us estimate the corresponding matrix of transition probabilities over two steps $P_{\mu_2}(i,j)$. We denote in the elements of the matrix $\mu_i = \mu(\{i\})$ for brevity, and also $\xi(\mu, \gamma)  = (\mu_{1} (\gamma \mu_{1} +$ $+ 0.001) + 0.499 \mu_{2} + 0.499 \mu_{3} + 0.001 \mu_{4})$ and $\zeta(\mu, \gamma, x) = (\gamma \mu_{1}+ $ $+0.001) \left(\gamma \xi(\mu, \gamma) + x\right)$ and $\delta(\gamma, \mu, x, y) = y\gamma\mu_1 + \zeta(\gamma, \mu, x)$, then the transition kernel over two steps has the form
{\small
\begin{align*}
P_{\mu_2}(i,j) =&
\left(\begin{matrix}
\delta(\gamma, \mu, 0.001, -0.001) + 0.249999 &
\delta(\gamma, \mu, 0.001, -0.001) + 0.249999 &
\dots\\
0.499 \gamma \xi(\mu, \gamma) + 0.25 &
0.499 \gamma \xi(\mu, \gamma) + 0.25 &
\dots\\
0.499 \gamma \xi(\mu, \gamma) + 0.25 &
0.499 \gamma \xi(\mu, \gamma) + 0.001996 &
\dots\\
0.001 \gamma \xi(\mu, \gamma) + 0.25 &
0.001 \gamma \xi(\mu, \gamma) + 0.498004 &
\dots\\
\end{matrix}\right.\\
&\left.\begin{matrix}
\dots &
\delta(\gamma, \mu, -0.499, -0.499) + 0.249501 &
\delta(\gamma, \mu, -0.499, -0.499) + 0.249501\\
\dots &
- 0.499 \gamma \xi(\mu, \gamma) + 0.25 &
- 0.499 \gamma \xi(\mu, \gamma) + 0.25\\
\dots &
- 0.499 \gamma \xi(\mu, \gamma) + 0.498004 &
- 0.499 \gamma \xi(\mu, \gamma) + 0.25\\
\dots &
- 0.001 \gamma \xi(\mu, \gamma) + 0.001996 &
- 0.001 \gamma \xi(\mu, \gamma) + 0.25
\end{matrix}\right)
\end{align*}
}
In this case $\alpha_2 = 0.503992$, $\lambda_2 \leqslant \gamma$. Therefore, we have $\lambda_2 < \alpha_2$, which leads to exponential convergence of the process.

\subsection{Example of a nonlinear Markov chain, $\alpha = 0, \lambda > \alpha$}
Let us show that this result can be used in cases when the one-step estimate is inapplicable, and, at the same time, violation of the \cite{Butkovsky:2014} conditions does not prevent some nonlinear Markov chains from exponential convergence. Consider the following discrete nonlinear Markov chain $X_n^\mu$ with state space $(E, \mathcal{E}) = \left(\{1, 2, 3, 4\}, 2^{\{1, 2, 3, 4\}}\right)$, the initial distribution $\mu_0$ and the transition probability matrix $P_{\mu_0}(i, j)$, defined as follows:
\vspace{-0.5em}
\begin{equation*}
\mu_0 = \begin{pmatrix}\mu(\{1\}) & \mu(\{2\})& \mu(\{3\})& \mu(\{4\})\end{pmatrix},
\vspace{-0.25em}\end{equation*}
\vspace{-0.5em}
\begin{equation*}
P_{\mu_0}(i, j) =
\begin{pmatrix}
0 & \gamma  \mu(\{1\}) & 0.5 - \gamma  \mu(\{1\}) & 0.5\\
0.5 & 0.5 & 0 & 0 \\
0.5 & 0 & 0.5 & 0 \\
0 & 0.5 & 0 & 0.5
\end{pmatrix},
\vspace{-0.25em}\end{equation*}
where $0 < \gamma < 0.5$.

We can notice that for a given process the conditions \eqref{NonlinearDobrushinAnalog} and \eqref{ButkovskyLambda} do not guarantee convergence to an invariant measure, since  $\lambda > \alpha$, since $\alpha = 0$ and $\lambda = \gamma$.

However, if we consider the corresponding matrix of transition probabilities in two steps $Q_{\mu_0}(i,j)$,
{\small
\begin{align*}
Q_{\mu_0}(i,j) =
\begin{pmatrix}
0.25 & 0.5 \gamma \mu_{1} + 0.25 & - 0.5 \gamma \mu_{1} + 0.25 & 0.25\\
0.25 & 0.25 \gamma \left(\mu_{2} + \mu_{3}\right) + 0.25 & - 0.25 \gamma \left(\mu_{2} + \mu_{3}\right) + 0.25 & 0.25\\
0.25 & 0.25 \gamma \left(\mu_{2} + \mu_{3}\right) & - 0.25 \gamma \left(\mu_{2} + \mu_{3}\right) + 0.5 & 0.25\\
0.25 & 0.5 & 0 & 0.25
\end{pmatrix},
\end{align*}}
we may obtain a different result. We have $\lambda_2 < \alpha_2$, as $\lambda_2 = \gamma/2$, whie $\alpha_2$ reaches its minimum at the pair of states $\{3, 4\}$ with a value in interval $[0.5; 0.5 + 0.25 \gamma]$.
Thus, the proposed estimate can guarantee exponential convergence in some cases when the existing \cite{Butkovsky:2014} result does not work.

\section{Conclusion}
The article proposes an improved estimate for the rate of convergence of homogeneous nonlinear Markov chains with discrete-time by generalizing the existing results on convergence and obtaining an estimate in several steps. This estimate leads to better convergence and may even be applicable in cases when one-step estimate cannot guarantee any convergence. In the last section of the work, an example of a homogeneous nonlinear Markov chain with a finite state space and discrete-time is given, which illustrates this result. In addition, this example shows that failure to satisfy the convergence conditions proposed in \cite{Butkovsky:2014} does not prevent the existence of a unique invariant measure and exponential convergence for homogeneous nonlinear Markov chains in discrete time.

\end{document}